\documentclass[12pt]{amsart}

\usepackage{amsmath,xspace,amssymb,mathrsfs}
\usepackage{color}

\input xy
\xyoption{all}
\xyoption{2cell}
\UseAllTwocells
\CompileMatrices

\newcommand{\Spec}{\operatorname{Spec}}
\renewcommand{\phi}{\varphi}

\newcommand{\Ker}{\operatorname{Ker}}

\newcommand{\Ima}{\operatorname{Im}}

\newtheorem{proposition}{Proposition}[section]
\newtheorem{lemma}[proposition]{Lemma} 
\newtheorem{corollary}[proposition]{Corollary}
\newtheorem{theorem}[proposition]{Theorem}
\newtheorem{bk}[proposition]{}

\newtheorem{prop-def}[proposition]{Proposition and definition}

\theoremstyle{definition}

\begin{document}

\title[Closed subschemes of a scheme]{On the closed subschemes of a scheme}

\author[A. Tarizadeh]{Abolfazl Tarizadeh}
\address{ Department of Mathematics, Faculty of Basic Sciences, University of Maragheh \\
P. O. Box 55136-553, Maragheh, Iran.
 }
\email{ebulfez1978@gmail.com}


\begin{abstract} In this paper, we obtain some new results on closed subschemes. Specially, we define natural addition and multiplication on the closed subschemes of a scheme. It is shown that ``the multiplication'' precisely coincides with the well known notion of ``the scheme-theoretic intersection''.  Dually, ``the addition'' coincides with ``the scheme-theoretic union''. It is also proved that these structures naturally provide contravariant functors from the category of schemes to the category of commutative monoids. \\
\end{abstract}

\maketitle

\section{Introduction}

It is a fundamental fact that the category of schemes admits all finite limits. Specially, it admits pullbacks. But this category is not well-behaved with finite colimits. For instance, pushouts do not necessarily exist in this category. Investigating the existence of pushouts in the category of schemes has been one of the main streams of algebraic geometry over the more recent years, see e.g. \cite{Ferrand}, \cite{Li}, \cite{Schwede}, \cite[Tag 0ECH]{Johan} and \cite{Temkin}. The most important result of theses articles
is that the pushout of any two closed immersions of schemes with a fixed source exists in the category of schemes. In the present paper, we use the whole strength of this result in order to define addition on the closed subschemes of a scheme. Dually, pullbacks of schemes allows us to define multiplication on the closed subschemes.\\

The collection of open subschemes of a scheme $X$ is in bijection to the set of open subsets of $X$, but a similar result does not hold for closed subschemes. As a specific example, consider the closed subset $Y=\{p\mathbb{Z}\}$ of $\Spec\mathbb{Z}$ where $p$ is a prime number. Then for each $n\geq1$ the canonical ring map $\mathbb{Z}\rightarrow\mathbb{Z}/p^{n}\mathbb{Z}$ induces a closed subscheme structure over $Y$. Thus $Y$ admits infinitely many pairwise distinct closed subscheme structures. In spite of this, the class of all closed subschemes of a scheme $X$ still forms a ``set", see Theorem \ref{Theorem 2}. In fact, if $X$ is an affine scheme then by Theorem \ref{Theorem 1}, the class of all closed subschemes of $X$ is in bijection to the set of ideals of $\mathscr{O}_{X}(X)$. The general case is a little complicated. We shall prove that a similar embedding hold in the general case. More precisely, we show that the class of all closed subschemes of a scheme $X$ can be canonically embedded into the set $\prod\limits_{U\in(U_{i})}S_{U}$ where $(U_{i})$ is an affine open covering of $X$ and each $S_{U}$ is the set of ideals of $\mathscr{O}_{X}(U)$. This result allows us to define addition (scheme-theoretic union) and multiplication (scheme-theoretic intersection) on the closed subschemes of a scheme. Then we show that these structures naturally provide contravariant functors from the category of schemes to the category of commutative monoids. \\

\section{Preliminaries}

In this section for convenience of the reader we briefly recall some material from the basic algebraic geometry which is needed in the next section. \\

\begin{bk} If $\phi:X\rightarrow Y$ and $\psi:Y\rightarrow Z$ are two morphisms of ringed spaces then $(\psi\circ\phi)^{\sharp}=\psi_{\ast}(\phi^{\sharp})\circ\psi^{\sharp}$ and $(\psi\circ\phi)_{x}^{\sharp}=\phi_{x}^{\sharp}\circ\psi_{\phi(x)}^{\sharp}$ for all $x\in X$. \\
\end{bk}

\begin{lemma}\label{Lemma 1} Let $f,g:Z\rightarrow Y$ be two morphisms of locally ringed spaces with $Y$ an affine scheme. If $f_{Y}^{\sharp}=g_{Y}^{\sharp}$ then $f=g$. \\
\end{lemma}

\begin{bk} Let $\phi:X\rightarrow Y$ be a morphism of ringed spaces, $V\subseteq Y$ and $U\subseteq X$ are opens such that $U\subseteq\phi^{-1}(V)$. Then $\phi$ induces a unique morphism of ringed spaces $\psi:U\rightarrow V$ such that the following diagram is commutative: $$\xymatrix{
U\ar[r]^{\psi} \ar[d]^{i} &V
\ar[d]^{j} \\X\ar[r]^{\phi} &Y}$$ where $i$ and $j$ are the inclusion morphisms. The morphism $\psi$ is often denoted by $\phi|_{U}$. If moreover $\phi$ is a morphism of locally ringed spaces then $\psi$ is as well. \\
\end{bk}

\begin{bk}\label{Blank 1} Let $f,g:Z\rightarrow Y$ be two morphisms of ringed spaces such that the maps $f$ and $g$ between the underlying spaces are the same. If there exists an open covering $(V_{i})$ of $Y$ such that $f|_{f^{-1}(V_{i})}=g|_{f^{-1}(V_{i})}:f^{-1}(V_{i})\rightarrow V_{i}$ for all $i$ then $f^{\sharp}=g^{\sharp}$ and so $f=g$. \\
\end{bk}

By a closed immersion we mean a morphism of ringed spaces $\phi:X\rightarrow Y$ such that the map $\phi$ between the underlying spaces is injective and closed map and for each $x\in X$ the ring map $\phi_{x}^{\sharp}:\mathscr{O}_{Y,\phi(x)}\rightarrow\mathscr{O}_{X,x}$ is surjective. If $\phi$ is a closed immersion and $V$ an open of $Y$ then the induced morphism $\phi|_{U}:U=\phi^{-1}(V)\rightarrow V$ is a closed immersion. A closed subscheme of a scheme $X$ is a pair $(Z,j)$ where $Z$ is a scheme and $j:Z\rightarrow X$ is a closed immersion such that the map $j$ between the underlying spaces is the canonical injection. It is simply denoted by $Z$ if there is no confusion on $j$. \\

\begin{theorem}\label{Theorem 1} Let $\phi:Y\rightarrow X$ be a closed immersion of schemes with $X$ an affine scheme. Then $Y$ is an affine scheme and there exists an isomorphism of schemes $g:Y\rightarrow\Spec A/I$ such that the following diagram is commutative:  $$\xymatrix{
Y\ar[r]^{\phi} \ar[d]^{g}&X\ar[d]^{\pi}\\\Spec A/I\ar[r]^{\:\:\:\:\:\eta}&\Spec A}$$ where $A=\mathscr{O}_{X}(X)$, $I=\Ker\phi_{X}^{\sharp}$, the morphism $\eta$ is induced by the canonical ring map $A\rightarrow A/I$ and $\pi$ is the canonical isomorphism. Moreover such an ideal $I$ is unique. \\
\end{theorem}

\begin{theorem}\label{Proposition 2} Closed immersions of schemes are stable under base change. \\
\end{theorem}

\begin{lemma}\label{Lemma 3} $($Glueing morphisms$)$ Let $Y$ and $Z$ be two ringed spaces and $(V_{i})$ an open covering of $Y$. Assume that for each $i$ a morphism $\theta_{i}:V_{i}\rightarrow Z$ of ringed spaces is given such that for each pair of indices $(i,j)$, $\theta_{i}|_{V_{i}\cap V_{j}}=\theta_{j}|_{V_{i}\cap V_{j}}$. Then there exists a unique morphism of ringed spaces $\theta:Y\rightarrow Z$ such that $\theta\circ j_{V_{i}}=\theta_{i}$ for all $i$ where $j_{V_{i}}:V_{i}\rightarrow Y$ is the inclusion morphism. If moreover all of the $\theta_{i}$ are morphisms of locally ringed spaces then $\theta$ is as well.\\
\end{lemma}

\begin{bk} Fix a scheme $X$. Let $\phi:Y\rightarrow X$ and $\psi:Z\rightarrow X$ be two closed immersions of schemes. We say that $\phi\sim\psi$ if there exists an isomorphism of schemes $\theta:Y\rightarrow Z$ such that $\phi=\psi\circ\theta$. We shall denote by $M(X)$ the collection of isomorphism classes of closed immersions of schemes with target $X$. \\
\end{bk}

\begin{proposition}\label{Corollary 1} The class of all closed subschemes of a  scheme $X$ is in bijection to $M(X)$.  \\
\end{proposition}

{\bf Proof.} Let $\phi:Y\rightarrow X$ be a closed immersion of schemes. Then $[\phi]$ contains a unique element $j:Z\rightarrow X$ such that $(Z,j)$ is a closed subscheme of $X$. Because let $Z:=\Ima\phi$ and $\mathscr{O}_{Z}:=f_{\ast}\mathscr{O}_{Y}$ where $f:Y\rightarrow Z$ is induced by $\phi$ and $f^{\sharp}:\mathscr{O}_{Z}\rightarrow f_{\ast}\mathscr{O}_{Y}$ is the identity morphism. If $U\subseteq X$
is an open then $j_{U}^{\sharp}:=\phi_{U}^{\sharp}$. Then clearly $f$ is an isomorphism of schemes, $j$ is a closed immersion and $\phi=j\circ f$. If $\psi:Y'\rightarrow X$ is a second closed imersion of schemes which induces the closed subscheme $(Z',j')$ as above. Then $(Z,j)=(Z',j')$ if and only if $[\phi]=[\psi]$. Therefore the map $[\phi]\rightsquigarrow(Z,j)$ is a bijection between $M(X)$ and the class of all closed subschemes of $X$. $\Box$ \\

\begin{theorem}\label{Theorem 4} Let $\phi:Y\rightarrow X$ and $\psi:Z\rightarrow X$ be two closed immersions of schemes. Then there exists a closed immersion $\theta:T\rightarrow X$ and morphisms $\theta_{1}:Y\rightarrow T$ and $\theta_{2}:Z\rightarrow T$ such that $\phi=\theta\circ\theta_{1}$, $\psi=\theta\circ\theta_{2}$ and the following: $$\xymatrix{
Y\times_{X}Z\ar[r]^{\:\:\:\:\:\:\phi'}\ar[d]^{\psi'} &Z\ar[d]^{\theta_{2}} \\Y\ar[r]^{\theta_{1}}&T}$$ is a pushout diagram where $\phi'$ and $\psi'$ are the canonical projections. If moreover $U\subseteq X$ is an affine open then the following diagram: $$\xymatrix{
\phi^{-1}(U)\times_{U}\psi^{-1}(U)\ar[r]^{}
\ar[d]&\psi^{-1}(U)\ar[d]^{} \\\phi^{-1}(U)\ar[r]^{}&\theta^{-1}(U)}$$ is a pushout diagram. \\
\end{theorem}

{\bf Proof.} To see its proof please consider \cite[Tag 0ECH]{Johan} specially \cite[Tag 0E25]{Johan} and \cite[Tag 0B7M]{Johan}. $\Box$ \\

\section{Main results}

\begin{bk} Fix a scheme $X$ and let $\mathscr{C}=(U_{i})$ be an affine open covering of $X$. \\
\end{bk}

\begin{lemma}\label{Lemma 2} Every monomorphism in the category of affine schemes is a monomorphism in the category of schemes. \\
\end{lemma}

{\bf Proof.} Let $\phi:Y\rightarrow X$ be a monomorphism in the category of affine schemes and $f,g:Z\rightarrow Y$ two morphisms of schemes such that $\phi\circ f=\phi\circ g$. We have $(\phi\circ f)_{X}^{\sharp}=(\phi\circ g)_{X}^{\sharp}$. It follows that $f_{Y}^{\sharp}\circ\phi_{X}^{\sharp}=g_{Y}^{\sharp}\circ\phi_{X}^{\sharp}$. But $\phi_{X}^{\sharp}$ is an epimorphism of rings because the category of affine schemes is anti-equivalent to the category of commutative rings. Thus
$f_{Y}^{\sharp}=g_{Y}^{\sharp}$. Therefore by Lemma \ref{Lemma 1}, $f=g$. $\Box$ \\

\begin{proposition}\label{Proposition 1} Every closed immersion of schemes is a monomorphism. \\
\end{proposition}

{\bf Proof.} Let $\phi:Y\rightarrow X$ be a closed immersion of schemes and
$f,g:Z\rightarrow Y$ two morphisms of schemes such that $\phi\circ f=\phi\circ g$. Then the maps $f,g$ between the underlying spaces are the same. To prove $f^{\sharp}=g^{\sharp}$ we act as follows. If $U$ is an affine open of $X$ then by Theorem \ref{Theorem 1}, $V=\phi^{-1}(U)$ is an affine open and the induced morphism $\phi|_{V}:V\rightarrow U$ is a monomorphism in the category of affine schemes. Thus by Lemma \ref{Lemma 2}, $\phi|_{V}$ is a monomorphism in the category of schemes and so $i_{U}\circ\phi|_{V}:V\rightarrow X$ is a monomorphism where $i_{U}:U\rightarrow X$ is the inclusion morphism. We have $i_{U}\circ\phi|_{V}\circ f|_{W}=i_{U}\circ\phi|_{V}\circ g|_{W}$ where $f|_{W},g|_{W}:W=f^{-1}(V)\rightarrow V$ are induced by $f$ and $g$, respectively. Thus $f|_{W}=g|_{W}$. So by \ref{Blank 1}, $f=g$. $\Box$ \\

\begin{theorem}\label{Theorem 2} The function $\eta:M(X)\rightarrow\prod\limits_{U\in\mathscr{C}}S_{U}$ given by  $[\phi]\rightsquigarrow(\Ker\phi_{U}^{\sharp})$ is injective. \\
\end{theorem}

{\bf Proof.} Let $\phi:Y\rightarrow X$ be a closed immersion. If $U\subseteq X$ is an affine open then by Theorem \ref{Theorem 1}, there exists a unique ideal $I=\Ker\phi_{U}^{\sharp}\subseteq A=\mathscr{O}_{X}(U)$ and an isomorphism of schemes $g_{U}:\phi^{-1}(U)\rightarrow\Spec A/I$ such that the following diagram is commutative: $$\xymatrix{
\phi^{-1}(U)\ar[r]^{\:\:\:\:\:\:\phi|_{\phi^{-1}(U)}} \ar[d]_{g_{U}}&U\ar[d]^{\pi}\\\Spec A/I\ar[r]^{\:\:\:\:\:\eta}&\Spec A}$$ where $\eta$ is induced by the canonical ring map $A\rightarrow A/I$. Assume that there is another closed immersion of schemes $\psi:Z\rightarrow X$ in which for each affine open $U\subseteq X$ similarly above there exists an isomorphism $h_{U}:\psi^{-1}(U)\rightarrow\Spec A/I$ such that the following diagram is commutative: $$\xymatrix{
\psi^{-1}(U)\ar[r]^{\:\:\:\:\:\:\psi|_{\psi^{-1}(U)}} \ar[d]_{h_{U}}&U\ar[d]^{\pi}\\\Spec A/I\ar[r]^{\:\:\:\:\:\eta}&\Spec A.}$$ We shall find an isomorphism $\theta:Y\rightarrow Z$ such that $\phi=\psi\circ\theta$. For each affine open $U\subseteq X$ take $\theta_{\phi^{-1}(U)}:=i_{\psi^{-1}(U)}\circ h^{-1}_{U}\circ g_{U}:\phi^{-1}(U)\rightarrow Z$ where $i_{\psi^{-1}(U)}:\psi^{-1}(U)\rightarrow Z$ is the inclusion morphism. If $U'\subseteq X$ is another affine open then we have $\psi\circ\big(\theta_{\phi^{-1}(U)}|_{V}\big)=
\phi\circ i_{V}=
\psi\circ\big(\theta_{\phi^{-1}(U')}|_{V}\big)$ where $V:=\phi^{-1}(U)\cap\phi^{-1}(U')$ and $i_{V}:V\rightarrow
 Y$ is the inclusion morphism. Thus by Proposition \ref{Proposition 1},
$\theta_{\phi^{-1}(U)}|_{V}=
\theta_{\phi^{-1}(U')}|_{V}$. Therefore by Lemma \ref{Lemma 3}, there exists a (unique) morphism of schemes $\theta:Y\rightarrow Z$ such that for each affine open $U\subseteq X$, $\theta\circ i_{\phi^{-1}(U)}=\theta_{\phi^{-1}(U)}$. Clearly the maps $\phi$ and $\psi\circ\theta$ between the underlying spaces are the same. It follows that $\theta$ is an isomorphism. Moreover $i_{U}\circ\big((\psi\circ\theta)|_{\phi^{-1}(U)}\big)=
i_{U}\circ(\phi|_{\phi^{-1}(U)})$. It follows that $(\psi\circ\theta)|_{\phi^{-1}(U)}=
\phi|_{\phi^{-1}(U)}$ because each open immersion is a monomorphism. Thus by \ref{Blank 1}, $\phi=\psi\circ\theta$. Hence $\eta$ is injective. $\Box$ \\

\begin{bk} If $I:X\rightarrow X$ is the identity morphism then $\eta$ maps $[I]$ to the sequence of zero ideals. Moreover $\eta$ maps the isomorphism class of the canonical morphism $\emptyset\rightarrow X$ to the sequence of unit ideals $\big(\mathscr{O}_{X}(U)\big)$ where $\emptyset$ is the empty scheme. \\
\end{bk}

\begin{bk} Let $\phi:Y\rightarrow X$ and $\psi:Z\rightarrow X$ be two closed immersions of schemes and consider the following pullback diagram:  $$\xymatrix{
Y\times_{X}Z\ar[r]^{\:\:\:\:\:\:\phi'}\ar[d]^{\psi'} &Z\ar[d]^{\psi} \\Y\ar[r]^{\phi}&X.}$$ By Theorem \ref{Proposition 2}, $\phi'$ is a closed immersion and so $\psi\circ\phi'=\phi\circ\psi'$ is a closed immersion. Thus we may define an operation over $M(X)$ as $[\phi].[\psi]=[\phi\circ\psi']$. It is called the scheme-theoretic intersection of $[\phi]$ and $[\psi]$.\\
\end{bk}

\begin{bk}\label{Blank 202} We define the operation $+$ on $M(X)$ as $[\phi]+[\psi]=[\theta]$, for $\theta$ see Theorem \ref{Theorem 4}. It is called the scheme-theoretic union of $[\phi]$ and $[\psi]$. \\
\end{bk}

\begin{theorem}\label{Theorem 3} If $\phi:Y\rightarrow X$ and $\psi:Z\rightarrow X$ are two closed immersions of schemes then the following hold.\\
$\mathbf{(i)}$ $\eta$ maps $[\phi].[\psi]$ into the sequence $(\Ker\phi^{\sharp}_{U}+\Ker\psi^{\sharp}_{U})$. \\
$\mathbf{(i)}$ $\eta$ maps $[\phi]+[\psi]$ into the sequence $(\Ker\phi^{\sharp}_{U}\cap\Ker\psi^{\sharp}_{U})$. \\
\end{theorem}

{\bf Proof.} The following: $$\xymatrix{
W=(\phi\circ\psi')^{-1}(U)\ar[r]^{}\ar[d]^{} &\psi^{-1}(U)\ar[d]^{} \\\phi^{-1}(U)\ar[r]^{}&U}$$ is a pullback diagram in the category of affine schemes because $W=\phi^{-1}(U)\times_{U}\psi^{-1}(U)$. By applying the global sections functor to the above diagram then we get the following pushout diagram in the category of commutative rings: $$\xymatrix{
A\ar[r]^{\psi_{U}^{\sharp}}\ar[d]^{\phi_{U}^{\sharp}} &C\ar[d]^{g} \\B\ar[r]^{f}&D}$$ where $A=\mathscr{O}_{X}(U)$, $B=\mathscr{O}_{Y}\big(\phi^{-1}(U)\big)$,
$C=\mathscr{O}_{Z}\big(\psi^{-1}(U)\big)$, $D=\mathscr{O}_{S}(W)$, $S=Y\times_{X}Z$, $f=\psi_{\phi^{-1}(U)}^{_{'}\sharp}$ and $g=\phi_{\psi^{-1}(U)}^{_{'}\sharp}$. By Theorem \ref{Theorem 1}, there exist isomorphisms of rings $\lambda:B\rightarrow A/I$ and $\mu:C\rightarrow A/J$ such that $\lambda\circ\phi_{U}^{\sharp}:A\rightarrow A/I$ and $\mu\circ\psi_{U}^{\sharp}:A\rightarrow A/J$ are the canonical maps where $I=\Ker\phi_{U}^{\sharp}$ and $J=\Ker\psi_{U}^{\sharp}$. Clearly the following: $$\xymatrix{
A\ar[r]\ar[d]^{} &A/J\ar[d]^{h_{2}}\\A/I\ar[r]^{h_{1}\:\:\:\:\:}&A/I+J}$$ is a pushout diagram in the category of commutative rings where $h_{1}$ and $h_{2}$ are the canonical maps. Thus there exists an isomorphism of rings $h:D\rightarrow A/I+J$ such that $h_{1}\circ\lambda=h\circ f$ and $h_{2}\circ\mu=h\circ g$. It follows that $\Ker(\phi\circ\psi')_{U}^{\sharp}=I+J$. This completes the proof of $(i)$. By Theorem \ref{Theorem 4}, the following: $$\xymatrix{
R\ar[r]^{g'}\ar[d]^{f'} &C\ar[d]^{g} \\B\ar[r]^{f}&D}$$ is a pullback diagram in the category of commutative rings where $R=\mathscr{O}_{T}\big(\theta^{-1}(U)\big)$, $f'=(\theta^{\sharp}_{1})_{\theta^{-1}(U)}$ and $g'=(\theta^{\sharp}_{2})_{\theta^{-1}(U)}$. Clearly the following: $$\xymatrix{
A/I\cap J\ar[r]^{\pi_{2}}\ar[d]^{\pi_{1}} &A/J\ar[d]^{h_{2}}\\A/I\ar[r]^{h_{1}\:\:\:\:\:}&A/I+J}$$ is a pullback diagram in the category of commutative rings where $\pi_{1}$ and $\pi_{2}$ are the canonical maps. Thus there exists an isomorphism of rings $h':R\rightarrow A/I\cap J$ such that $\lambda\circ f'=\pi_{1}\circ h'$ and $\mu\circ g'=\pi_{2}\circ h'$. It follows that $h'\circ\theta_{U}^{\sharp}:A\rightarrow A/I\cap J$ is the canonical map and so $\Ker\theta_{U}^{\sharp}=I\cap J$. $\Box$ \\

\begin{corollary}\label{Corollary 2} The set $M(X)$ with the operation $.$ forms a commutative monoid such that $[\phi].[\phi]=[\phi]$ for all $[\phi]\in M(X)$. \\
\end{corollary}

{\bf Proof.} It is an immediate consequence of Theorems \ref{Theorem 2} and \ref{Theorem 3}. $\Box$ \\

\begin{bk} Let $f:X\rightarrow Y$ be a morphism of schemes. Let $[\phi]\in M(Y)$ where $\phi:Z\rightarrow Y$ is a closed immersion. Consider the following pullback diagram: $$\xymatrix{
Z\times_{Y}X\ar[r]^{\:\:\:\:\:\:\phi'}\ar[d]^{p} &X\ar[d]^{f} \\Z\ar[r]^{\phi}&Y.}$$
By Theorem \ref{Proposition 2}, $\phi'$ is a closed immersion. Thus we may define $M(f):M(Y)\rightarrow M(X)$ as $[\phi]\rightsquigarrow[\phi']$. \\
\end{bk}

\begin{theorem}\label{Theorem 5} The assignments $X\rightsquigarrow M(X)$ and $f\rightsquigarrow M(f)$ define a contravariant functor from the category of schemes to the category of commutative monoids. \\
\end{theorem}

{\bf Proof.} First we show that $M(f)$ is a morphism of monoids. Clearly it preserves identities. Take a second element $[\psi]\in M(Y)$ where $\psi:T\rightarrow Y$ is a closed immersion and consider the following pullback diagrams:  $$\xymatrix{
T\times_{Y}X\ar[r]^{\:\:\:\:\:\:\psi'}\ar[d]^{p'} &X\ar[d]^{f} \\T\ar[r]^{\psi}&Y}$$ and  $$\xymatrix{
Z\times_{Y}T\ar[r]^{\:\:\:\:\:\:h}\ar[d]^{g} &T\ar[d]^{\psi} \\Z\ar[r]^{\phi}&Y.}$$ By definition, $[\phi]\ast[\psi]=[\phi\circ g]$. Set $A=Z\times_{Y}T$, $B=Z\times_{Y}X$ and $C=T\times_{Y}X$. Then consider the following pullback diagrams: $$\xymatrix{
A\times_{Y}X\ar[r]^{\:\:\:\:\:\:\:\:\:(\phi\circ g)'}\ar[d]^{p''} &X\ar[d]^{f} \\A\ar[r]^{\phi\circ g}&Y}$$ and $$\xymatrix{
B\times_{X}C\ar[r]^{\:\:\:\:\:\:\:\:\:q'''}\ar[d]^{p'''} &C\ar[d]^{\psi'} \\B\ar[r]^{\phi'}&X.}$$ We have to show that $[\phi'\circ p''']=[(\phi\circ g)']$. There exists a unique morphism $\mu:B\times_{X}C\rightarrow A$ such that $p\circ p'''=g\circ\mu$ and $p'\circ q'''=h\circ\mu$. It follows that there exists a unique morphism $\delta:B\times_{X}C\rightarrow A\times_{Y}X$
such that $\mu=p''\circ\delta$ and $\phi'\circ p'''=(\phi\circ g)'\circ\delta$. We must prove that $\delta$ is an isomorphism. There exists a unique morphism $\alpha:A\times_{Y}X\rightarrow B$ such that $g\circ p''=p\circ\alpha$ and $(\phi\circ g)'=\phi'\circ\alpha$. Also there exists a unique morphism $\beta:A\times_{Y}X\rightarrow C$ such that $h\circ p''=p'\circ\beta$ and $(\phi\circ g)'=\psi'\circ\beta$. Thus there exists a unique morphism $\lambda:A\times_{Y}X\rightarrow B\times_{X}C$ such that $\alpha=p'''\circ\lambda$ and $\beta=q'''\circ\lambda$. We have $\phi'\circ\alpha\circ\delta=\phi'\circ p'''$ and $\psi'\circ\beta\circ\delta=\psi'\circ q'''$. It follows that $\alpha\circ\delta=p'''$ and $\beta\circ\delta=q'''$ because by Proposition \ref{Proposition 1}, $\phi'$ and $\psi'$ are monomorphisms. Therefore $p'''\circ(\lambda\circ\delta)=\alpha\circ\delta=p'''$ and $q'''\circ(\lambda\circ\delta)=\beta\circ\delta=q'''$. It follows that $\lambda\circ\delta$ is the identity morphism. Clearly $(\phi\circ g)'\circ(\eta\circ\lambda)=(\phi\circ g)'$. We also have $\phi\circ g\circ(p''\circ\eta\circ\lambda)=\phi\circ g\circ p''$. It follows that $p''\circ(\eta\circ\lambda)=p''$ because by Theorem \ref{Proposition 2} and Proposition \ref{Proposition 1}, $\phi\circ g$ is a monomorphism. Thus $\delta\circ\lambda$ is the identity morphism. Therefore $M(f)$ is a morphism of monoids. If $I_{X}$ is the identity morphism of $X$ then clearly $M(I_{X})=I_{M(X)}$. It remains to show that if $f':Y\rightarrow Z$ is a second morphism of schemes then $M(f'\circ f)=M(f)\circ M(f')$. Let $\phi:T\rightarrow Z$ be a closed immersion of schemes and consider the following pullback diagrams: $$\xymatrix{
T\times_{Z}X\ar[r]^{\:\:\:\:\:\:\:\:\phi'}\ar[d]^{p} &X\ar[d]^{f'\circ f} \\T\ar[r]^{\phi}&Z,}$$  $$\xymatrix{
T\times_{Z}Y\ar[r]^{\:\:\:\:\:\:\:\:\phi''}\ar[d]^{p'} &Y\ar[d]^{f'} \\T\ar[r]^{\phi}&Z}$$ and  $$\xymatrix{
A\times_{Y}X\ar[r]^{\:\:\:\:\:\:\:\:\phi'''}\ar[d]^{p''} &X\ar[d]^{f} \\A\ar[r]^{\phi''}&Y}$$ where $A=T\times_{Z}Y$. There exists a canonical isomorphism $\lambda:A\times_{Y}X\rightarrow T\times_{Z}X$ such that $\phi'''=\phi'\circ\lambda$. Therefore $[\phi']=[\phi''']$. $\Box$ \\

Similarly above (i.e. Corollary \ref{Corollary 2}), the set $M(X)$ with the operation $+$, defined in \ref{Blank 202}, forms a commutative monoid such that $[\phi]+[\phi]=[\phi]$ for all $[\phi]\in M(X)$. We denote this monoid by $M^{+}(X)$. We have then the following result. \\

\begin{theorem} The assignments $X\rightsquigarrow M^{+}(X)$ and $f\rightsquigarrow M(f)$ define a contravariant functor from the category of schemes to the category of commutative monoids. \\
\end{theorem}

{\bf Proof.} It is proven exactly like Theorem \ref{Theorem 5}. $\Box$ \\

\end{document}